\magnification=1200
\input amstex
\documentstyle{amsppt}
\def\bl {\bold\Lambda}
\def\wt#1{\widetilde {#1}}

\def\SO {\Cal O}

\def\SE {\Cal E}
\def\SL {\Cal L}
\def\SF {\Cal F}

\def\SP {\Cal P}

\def\SR {\Cal R}
\def\SU {\Cal U}
\def\SW {\Cal W}
\def\SK {\Cal K}

\def\bset {\text{\bf Set}}
\topmatter
\title Degeneration of $SL(n)$-bundles on a reducible curve
\endtitle
\author Xiaotao Sun   \endauthor
\address Institute of Mathematics, Academia Sinica,
Beijing 100080, China \endaddress \email
xsun$\@$math08.math.ac.cn\endemail
\address Department of Mathematics, The University of 
Hong Kong, Pokfulam Road, Hong Kong\endaddress \email
xsun$\@$maths.hku.hk\endemail
\thanks The work supported by a grant of NFSC for outstanding
young researcher at contract number 10025103.
\endthanks
\endtopmatter
\document

\heading  Introduction \endheading

It is a classic idea in algebraic geometry to use degeneration
method. In particular, it achieved successes recently in the
studying of moduli spaces of vector bundles (See [Gi], [GL1],
[GL2], [NR] and [S1]). In connection of string theory, it needs
also to study the degeneration of moduli spaces of $G$-bundles for
any reductive algebraic group $G$ (See [F1],[F2]).

Let $\Cal X\to B$ be a proper flat family of curves of genus $g$
such that $\Cal X_b$ ($b\neq 0$) smooth and $\Cal X_0$ a
semistable curve. It was known that there exists a family $$\Cal
M_{\Cal X}(G)_0\to B_0=B\setminus\{0\}$$ of moduli spaces of
semistable $G$-bundles. The question becomes that for what
geometric objects on $\Cal X_0$ such that the moduli space of them
gives a compactification $\Cal M_{\Cal X}(G)\to B$ of $\Cal
M_{\Cal X}(G)_0\to B_0=B\setminus\{0\}$. We can consider $G$ as a
subgroup of some $GL(n)$, and think a semistable $G$-bundle as a
semistable vector bundle of rank $n$ with some additional
conditions. Thus we may think (not strictly) the moduli space of
semistable $G$-bundles as a subscheme of the moduli space of
semistable vector bundles of rank $n$. On the other hand, there is
a natural choice of geometric objects, the torsion free sheaves,
on $\Cal X_0$. The moduli space of semistable torsion free sheaves
gives a natural degeneration of moduli spaces of semistable vector
bundles on $\Cal X_b$ when $b$ goes to $0$. Then a possible
approach to the problem is finding the correct torsion free
$G$-sheaves on $\Cal X_0$. However, the problem remains almost
complete open except for special groups like $G=GL(n)$, $Sp(n)$
and $O(n)$ (See the introduction of [F1]). G. Faltings studied the
cases that $G=Sp(n)$ and $O(n)$ but left the case $G=SL(n)$ open
(See [F2]). In this paper, we will treat the case $G=SL(n)$ when
$\Cal X_0$ has two smooth irreducible components intersecting at
one node $x_0$.

Let $\SU_{\Cal X}\to B$ be the family of moduli spaces of
semistable torsion free sheaves and $S\SU_{\Cal X}^0\to B_0$ be
the family of moduli spaces of semistable vector bundles with
fixed determinant $\SL$. Let $f:S\SU_{\Cal X}\to B$ be the Zariski
closure of $S\SU_{\Cal X}^0\subset\SU_{\Cal X}$ in $\SU_{\Cal X}$.
Then the problem becomes to give a moduli interpretation of
$f^{-1}(0)$. Namely, to define a suitable moduli functor
$S\SU_{\Cal X_0}^{\natural}$ such that $f^{-1}(0)$ universally
corepresents $S\SU_{\Cal X_0}^{\natural}.$ It is obvious that the
above question can also be asked for moduli of bundles on higher
dimensional variety.  In the study of moduli spaces of bundles on
surfaces (See [GL1] and [GL2]),  Gieseker and Li have noted that
the closed points $[\SF]\in f^{-1}(0)$ satisfy the condition
$$det(\SF|_{\Cal X_0\setminus\{x_0\}})=\SL|_{\Cal
X_0\setminus\{x_0\}}.\tag{$*$}$$ In general, it may not be true
that a semistable sheaf $\SF$ satisfying condition $(*)$ has to be
a point of $f^{-1}(0)$. For example, when $\Cal X_0$ is
irreducible, the set of semistable sheaves satisfying the
condition $(*)$ will have bigger dimension than $f^{-1}(0)$.
However, in our case when $\Cal X_0$ has two smooth irreducible
components intersecting at one node, the points of $f^{-1}(0)$ are
precisely the semistable sheaves satisfying condition $(*)$ (See
Lemma 2.2). In fact, we defined a moduli functor $S\SU_{\Cal
X_0}^{\natural}$, which is represented by a closed subscheme
$S\SU_{\Cal X_0}\subset \SU_{\Cal X_0}$ of the moduli space of
semistable torsion free sheaves on $\Cal X_0$ (See Theorem 1.6).
Moreover, we proved that $S\SU_{\Cal X_0}$ is a reduced,
seminormal variety whose closed points are precisely the
$s$-equivalent classes of semistable sheaves on $\Cal X_0$
satisfying $(*)$. These was done in Section 1 (See Theorem 1.6).
In Section 2, we showed that the above moduli problem has good
specialzation (See [NS] for the notation), and the degeneration of
moduli spaces of semistable $SL(n)$-bundles is $S\SU_{\Cal X_0}$
when $b$ goes to $0$.

The paper comes from a conversation with Jun Li, who told me
the condition $(*)$ in his joint works with Gieseker. The emails
with D.S. Nagaraj and C.S. Seshadri concerning Proposition (4.1) of
[NS] were helpful for the observation of Lemma 1.4. I thank them 
very much.  

\heading\S1 Moduli space of semistable sheaves with fixed
determinant \endheading

Let $X$ be a projective curve of genus $g$ with two smooth
irreducible components $X_1$ and $X_2$ of genus $g_1$ and $g_2$,
which intersect at a node $x_0$ of $X$. Let
$X^0=X\smallsetminus\{x_0\}$ and $\SO_X(1)$ be a fixed ample line
bundle on $X$. Fix integers $r>0$, $d$ and $\chi=d+r(1-g)$.

For any sheaf $E$ on $X$,
$P(E,n)=\chi(E(n)):=dim\text{H}^0(E(n))-dim\text{H}^1(E(n))$ is
called Hilbert polynomial of $E$. If $E$ is a torsion free sheaf
on $X$, and let $r_i$ denote the rank of restriction of $E$ to
$X_i$ ($i=1,2$). The rank of $E$ is defined to be
$$rk(E):=\frac{1}{deg(\SO_X(1))}\cdot\lim_{n\to\infty}\frac{P(E,n)}{n}.$$
Let $c_i=deg(\SO_X(1)|_{X_i})$, $a_i=\frac{c_i}{c_1+c_2}$, then
$$P(E,n)=(c_1r_1+c_2r_2)n+\chi(E),\quad rk(E)=a_1r_1+a_2r_2.$$

\proclaim{Definition 1.1} A sheaf $E$ on $X$ is called semistable
(resp. stable) if, for any subsheaf $E_1\subset E$, one has
$$\chi(E_1)\le\,\,(resp.\,\,<)\,\,\,\frac{\chi(E)}{rk(E)}\cdot
rk(E).$$ \endproclaim

The moduli functor $\SU_X^{\natural}:(\Bbb C-schemes)\to (sets)$
was defined as $$\SU_X^{\natural}(S)=\left\{\aligned
&\text{$\SO_S$-flat semistable sheaves $E$ on $X\times S$}\\
&\text{of rank $r$ and $\chi(E|_{X\times\{s\}})=\chi$ for any
$s\in S$}
\endaligned\right\}$$
and it is known that there exists a projective scheme $\SU_X$,
which universally corepresents the functor $\SU_X^{\natural}$.

For any integer $N$ and polynomial $P(n)=r(c_1+c_2)n+\chi$, let
$\SW=\SO_X(-N)^{P(N)}$ and $\SW_S=\SO_{X\times S}(-N)^{P(N)}$.
Recall that $Quot_{P(n)}(\SW)$ is Grothendieck quotient scheme,
whose $S$-valued points may be described as the set of quotients
$$\SW_S\twoheadrightarrow F$$ on $X\times S$, where $F$ is flat
over $S$ and its Hilbert polynomial is $P(n)$. Let
$\SR^{ss}\subset Quot_{P(n)}(\SW)$ (resp. $\SR^s$) be the open set
where the sheaf $F$ is semistable (resp. stable) and
$\SW_S\twoheadrightarrow F$ induces an isomorphism
$H^0(\SW_S(N))\cong H^0(F(N)).$ Since the set of semistable
sheaves with fixed Hilbert polynomial is bounded, we can assume
that $N$ is chosen large enough so that: every semistable sheaf
with Hilbert polynomial $P(n)$ appears as a point in $\SR^{ss}$.
The group $SL(P(N))$ acts on $Quot_{P(n)}(\SW)$ and thus on
$\SR^{ss}$. The moduli space $\SU_X$ was constructed as a good
quotient $$\SU_X=\SR^{ss}//SL(P(N)).$$ More precisely, the
following were known (See [Si] for the more general results)

\proclaim{Theorem 1.2} ([Se], [Si]) Let $\SU_X=\SR^{ss}//SL(P(N))$
be the good quotient. Then \roster
\item There exists a natural transformation
$\SU_X^{\natural}\to\SU_X$ such that $\SU_X$ universally
corepresents $\SU_X^{\natural}$.
\item $\SU_X$ is projective, and its geometric points are
$$\SU_X(\Bbb C)=\left\{\aligned &\text{$s$-equivalent classes of
semistable}\\ &\text{sheaves of rank $r$ and degree $d$}
\endaligned\right\}.$$
\item There is an open subset $\SU_X^s\subset \SU_X$, with inverse
image equal to $\SR^s$, whose points represent isomorphism classes
of stable sheaves. Locally in the {\'e}tale topology on $\SU_X^s$,
there exists a universal sheaf $\SE^{univ}$ such that if
$\SE\in\SU_X^{\natural}(S)$ whose fibres $\SE_s$ are stable, then
the pull-back of $\SE^{univ}$ via $S\to\SU^s_X$ is isomorphic to
$\SE$ after tensoring with the pull-back of a line bundle on
$S$.\endroster\endproclaim

Let $\SR_0\subset \SR^{ss}$ be the dense open set of locally free
sheaves. For each $F\in \SR_0$, let $F_1=F|_{X_1}$,
$F_2=F|_{X_2}$, we have $\chi(F_1)+\chi(F_2)=\chi +r$ and (by
semistability) $$a_1\chi\le\chi(F_1)\le a_1\chi +r,\quad
a_2\chi\le\chi(F_2)\le a_2\chi +r.\tag 1.1$$ Thus $\SR_0$ is the
disjoint union of
$\SR_0^{\chi_1,\chi_2}=\{F\in\SR_0|\,\,\chi(F|_{X_i})=\chi_i\}$
where $\chi_1+\chi_2=\chi+r$ satisfying $(1.1)$. Let
$\SR^{\chi_1,\chi_2}$ be the Zariski closure of
$\SR_0^{\chi_1,\chi_2}$ in $\SR^{ss}$, then
$$\SR^{ss}=\bigcup_{\chi_1,\chi_2}\SR^{\chi_1,\chi_2}$$ has at
most $r+1$ irreducible components. Let $\SU_{X_i}$ be the moduli
space of semistable vector bundles $F_i$ of rank $r$ and
$\chi(F_i)=\chi_i.$ It was known that for all possible choices of
$\chi_1$, $\chi_2$ satisfying $(1.1)$, $\SR^{\chi_1,\chi_2}$ is
not empty if $\SU_{X_1}$ and $\SU_{X_2}$ are not empty.
$$\SU_X=\bigcup_{\chi_1,\chi_2}\SU^{\chi_1,\chi_2}=
\bigcup_{\chi_1,\chi_2}\SR^{\chi_1,\chi_2}//SL(P(N))$$ has
$\sharp\{(\chi_1,\chi_2)\}$ irreducible components.

Let $\SL_0$ be a line bundle of degree $d$ on $X$, we define a
subfunctor $S\SU_X^{\natural}$ of $\SU_X^{\natural}$ by $$
S\SU_X^{\natural}(S)=\left\{\aligned &\text{$\SO_S$-flat
semistable sheaves $E$ on $X\times S$ of rank $r$
satisfying}\\&\text{ $det(E|_{X^0_S})=p_X^*(\SL_0)|_{X^0_S}$ and
$\chi(E|_{X\times\{s\}})=\chi$ for any $s\in S$}
\endaligned\right\}$$

We will prove that there exists a closed subscheme $S\SU_X\subset
\SU_X$ which is reduced, and universally corepresents the moduli
functor $S\SU_X^{\natural}$. To do that, we first prove that the
local deformation functor of it is pro-representable. Let $\bl$ be
the category of Artinian local $\Bbb C$-algebras. Morphisms in
$\bl$ are local homomorphisms of $\Bbb C$-algebras. For any $A\in
\bl$, we denote $X\times Spec(A)$ (resp. $X^0\times Spec(A)$) by
$X_A$ (resp. $X^0_A$). At any point $e_0=(\SW\twoheadrightarrow
E_0)$  of the $Quot$ scheme $Quot_{P(n)}(\SW)$, we have the local
deformation functor $$G(A)=\{\text{$A$-flat quotients
$(\SW_A\twoheadrightarrow E)$, with $\chi(E(n))=P(n)$}\}$$ such
that $G(\Bbb C)=\{e_0\}$. It is well known that
$\hat\SO_{Quot,e_0}$ pro-represents $G$. Let $\SL_0$ be a fixed
line bundle on $X$ such that $det(E_0|_{X^0})=\SL_0|_{X^0}$. We
can define a subfunctor of $G$ by
$$F(A)=\{(\SW_A\twoheadrightarrow E)\in G(A)\,|\,
det(E|_{X^0_A})=p_{X^0}^*\SL_0\},$$ where $p_{X^0}$ (resp. $p_X$,
$p_A$) denotes the projection to $X^0$ (resp. $X$, $Spec(A)$).

\proclaim{Proposition 1.3} The functor $F$ is pro-representable.
\endproclaim

\demo{Proof} Let $A_1\to A$ and $A_2\to A$ be morphisms in $\bl$,
and consider the map $$F(A_1\times_AA_2)\to
F(A_1)\times_{F(A)}F(A_2).\tag{1.2}$$ By Theorem 2.11 of [Sc], it
is enough to show that the map $(1.2)$ is bijective when $A_2\to
A$ is a small extension and the tangent space $t_F$ is a finite
dimensional vector space. But the bijectivity of $(1.2)$ implies
that $t_F$ has a vector space structure compatible with that of
$t_G$ (See Remark 2.13 of [Sc]). Thus $dim_{\Bbb C}(t_F)<\infty$
and we only need to check the bijectivity of (1.2).

Let $(\SW_{A_1}\twoheadrightarrow E_1)\in F(A_1),\quad(\SW_A
\twoheadrightarrow E)\in F(A)\quad (\SW_{A_2}\twoheadrightarrow
E_2)\in F(A_2)$ such that the restriction morphisms
$E_1@>{u_1}>>E$, $E_2@>{u_2}>>E$ induce isomorphisms
$$(\SW_{A_1}\twoheadrightarrow E_1)|_{X_A}\cong
(\SW_A\twoheadrightarrow E),\quad (\SW_{A_2}\twoheadrightarrow
E_2)|_{X_A}\cong (\SW_A\twoheadrightarrow E).$$ Let
$B=A_1\times_AA_2$, then, since $G$ is pro-representable, there
exists a unique
 $(\SW_B\twoheadrightarrow\SE)\in G(B)$ such that
$$(\SW_B\twoheadrightarrow
\SE)|_{X_{A_1}}\cong(\SW_{A_1}\twoheadrightarrow E_1),\quad
 (\SW_B\twoheadrightarrow\SE)|_{X_{A_2}}\cong(\SW_{A_2}\twoheadrightarrow E_2),$$
 where $\cong$ means equality as points of $Quot$ scheme. In
 particular, there exist morphism
 $\SE@>{q_1}>>E_1$, $\SE@>{q_2}>>E_2$ such that $u_1\cdot
 q_1=u_2\cdot q_2$ and $q_1$, $q_2$ induce isomorphisms
 $\SE|_{X_{A_1}}\cong E_1$, $\SE|_{X_{A_2}}\cong E_2$.
 By restricting everything to $X^0_B$ and taking wedge product, we
 have (for $i=1,2$)
 $$det(\SE^0)@>{det(q_i^0)}>>det(E_i^0),\quad
 det(E^0_i)@>{det(u_i^0)}>>det(E^0)$$
 satisfying $det(u^0_1)\cdot det(q_1^0)=det(u_2^0)\cdot
 det(q^0_2)$. Thus, by Corollary 3.6 of [Sc],
 $$det(\SE^0)\cong det(E^0_1)\times_{det(E^0)}det(E^0_2).$$

 To prove that $det(\SE^0)\cong p^*_{X^0}(\SL)$, we imitate the
 arguments of uniqueness in the proof of Proposition 3.2. Since
 $p^*_{X^0}(\SL)|_{X^0_{A_i}}\cong det(E^0_i)$, we have morphisms
 $$p^*_{X^0}(\SL)@>p_1>>det(E^0_1),\quad
 p^*_{X^0}(\SL)@>p_2>>det(E_2^0)$$
 which induce the isomorphisms and thus a commutative diagram
 $$\CD
det(E^0_2) @>{det(u^0_2)}>>  det(E^0)\\
 @Ap_2AA                                  \\
p^*_{X^0}(\SL)  @.              @A\theta AA         \\
 @Vp_1VV                              \\
det(E^0_1) @>det(u_1^0)>>     det(E^0)
\\
\endCD$$
where $\theta$ is an automorphism of $det(E^0)$. If there exists a
morphism $$det(E^0_2)@>\theta_2>>det(E^0_2)$$ such that $$\CD
det(E^0_2) @>\theta_2>>  det(E^0_2)  \\ @V det(u^0_2)VV @V
det(u^0_2)VV \\ det(E^0) @>\theta>> det(E^0) \\
\endCD$$
is commutative, then $\theta_2$ has to be an isomorphism (since
$A_2\to A$ is a small extension) by Lemma 3.3 of [Sc]. Thus we can
modify the morphism $p^*_{X^0}(\SL)@>p_2>>det(E_2^0)$ to $$\wt
p_2: p^*_{X^0}(\SL)@>p_2>>det(E_2^0)@>\theta_2^{-1}>>det(E_2^0),$$
so that $det(u^0_1)\cdot p_1=det(u_2^0)\cdot\wt p_2$. Thus
$det(\SE^0)\cong p^*_{X^0}(\SL)$ and we are done if the lift of
$\theta$ is always possible. But this is equivalent to the
surjectivity of the canonical map $$H^0(X^0_{A_2},
\SO_{X^0_{A_2}})\to H^0(X^0_A, \SO_{X^0_A}),$$ which is true since
for any finite dimensional $\Bbb C$-algebra $A$, we have
$$H^0(X^0_A,\SO_{X^0_A})=H^0(X^0,\SO_{X^0})\otimes_{\Bbb C}A.$$
\enddemo

For any point $e_0=(\SW\twoheadrightarrow E_0)$, we observe that
$F$ has the same singularity with $G$ at $e_0$. To see it, we
define a functor $T:\bl\to \bset$ by $$T(A)=\{\text{Isomorphism
classes of $A$-flat torsion free $(\SO_{X,x_0}\otimes
A)$-modules}\}$$ such that $T(\Bbb C)=\{E_0\otimes\SO_{X,x_0}\}.$
There is a morphism of functors $\phi: G\to T$ defined by
$$\phi((\SW\twoheadrightarrow E_A))=E_A\otimes(\SO_{X,x_0}\otimes
A).$$ It is known that $\phi$ is formally smooth (See Theorem 4.1
of [F2], or Proposition (4.1) of [NS]). Here we remark that its
restriction to the subfunctor $F$ is also formally smooth.

\proclaim{Lemma 1.4} The morphism $\phi: F\to T$ is formally
smooth.\endproclaim

\demo{Proof} Let $B\to A$ be a small extension, one need to check
the surjectivity of $$F(B)\to F(A)\times_{T(A)}T(B).$$ For any
$(\SW_A\twoheadrightarrow E_A)\in F(A)$ and $N\in T(B)$ satisfying
$E_A\otimes(\SO_{X,x_0}\otimes A)\cong N\otimes_BA,$ we can find
an open cover of $X$ consists two affine sets $U_1$, $U_2$ such
that $x_0\in U_2\smallsetminus U_1$ and $U_1\otimes Spec(A)$,
$(U_2\smallsetminus\{x_0\})\otimes Spec(A)$ trivializing the
vector bundle $E^0_A:=E_A|_{X^0_A}$ (See Proposition (4.1) of
[NS]). Thus the gluing data of $E^0_A$ is a matrix $M\in
GL(\SO_X(U_1\cap U_2)\otimes A)$ and $E_A$ is obtained by gluing
$E^0_A$ and $N\otimes_BA$. Then the lift $E_B\in F(B)$ was
obtained by lifting the gluing data (See Proposition (4.1) of
[NS]). Since $SL(\SO_X(U_1\cap U_2)\otimes B)\to SL(\SO_X(U_1\cap
U_2)\otimes A)$ is surjective, it is clear that we can choose a
lift $\wt M\in GL(\SO_X(U_1\cap U_2)\otimes B)$ of $M$ such that
the resulting sheaf $E_B\in F(B).$ This proves the lemma.
\enddemo

\proclaim{Corollary 1.5} The functor $F$ is pro-represented by a
reduced semi-normal complete local $\Bbb C$-algebra.
\endproclaim

\demo{Proof} This follows Lemma 1.4 and the study of functor $T$
in [F2] and [Se].\enddemo

To construct the closed subscheme $S\SU_X$, which will universally
corepresent the moduli functor $S\SU_X^{\natural}$, we recall some
constructions in [S1] and [S2]. Let $\pi:\wt X\to X$ be the
normalization of $X$ and $\pi^{-1}(x_0)=\{x_1,x_2\}$, then $\wt X$
is a disjoint union of $X_1$ and $X_2$ (we will identify $x_1$,
$x_2$ with $x_0$ when we work on $X$). A GPB $(E,Q)$ of rank $r$
on $\wt X$ is a vector bundle $E$ of rank $r$ on $\wt X$ (its
restriction to $X_i$ is denoted by $E_i$), together with a
quotient $E_{x_1}\oplus E_{x_2}\to Q$ of dimension $r$. We have
constructed the moduli space $$\SP:=
\coprod_{\chi_1+\chi_2=\chi+r}\SP_{\chi_1,\chi_2}$$ of
$s$-equivalence classes of semistable GPB $(E,Q)$ on $\wt X$ of
rank $r$ and $\chi(E)=\chi +r$ (See [S2]), where
$\SP_{\chi_1,\chi_2}=\{(E,Q)\in\SP\,|\,\chi(E_i)=\chi_i\}$ and
$\chi_1$, $\chi_2$ satisfy $(1.1).$ There are also finite
morphisms (See [S2])
$$\phi_{\chi_1,\chi_2}:\SP_{\chi_1,\chi_2}\to\SU_X^{\chi_1,\chi_2}\subset
\SU_X$$ such that $\phi=\coprod\phi_{\chi_1,\chi_2}: \SP\to\SU_X$
is the normalization of $\SU_X$. The morphism $\phi$ is defined
such that $\phi([(E,Q)]):=F$ satisfies  the exact sequence  $$0\to
F\to\pi_*E\to\, _{x_0}Q\to 0.$$ A straightforward generalization
of Lemma 5.7 in [S1] shows that there exists a morphism
$$Det:\SP_{\chi_1,\chi_2}\to J_{X_1}^{\chi_1-r(1-g_1)}\times
J_{X_2}^{\chi_2-r(1-g_2)}$$ such that
$Det([(E,Q)])=(det(E_1),det(E_2))$. Let $L_i=\SL_0|_{X_i}$
($i=1,2$) and $$L_1(n_1):=L_1\otimes\SO_{X_1}(n_1x_1),\quad
L_2(n_2):=L_2\otimes\SO_{X_2}(n_2x_2),$$ where
$n_i=\chi_i-r(1-g_i)-deg(L_i)$. We have the closed subschemes
$$\SP_{\chi_1,\chi_2}^{n_1,n_2}:=Det^{-1}((L_1(n_1),L_2(n_2)),$$
where $n_1$, $n_2$ are determined uniquely by $\chi_1$, $\chi_2$
satisfying $(1.1)$ since $\SL_0$ was fixed. Thus we denote the
closed subsets $\phi_{\chi_1,\chi_2}
(\SP_{\chi_1,\chi_2}^{n_1,n_2})\subset \SU_X^{\chi_1,\chi_2}$ by
$S\SU_X^{\chi_1,\chi_2}$, and define $S\SU_X$ to be the closed
subset
$$\bigcup_{\chi_1+\chi_2=\chi+r}S\SU_X^{\chi_1,\chi_2}\subset\SU_X$$
with the reduced scheme structure. Then we have

\proclaim{Theorem 1.6} Let $S\SU_X$ be the closed subscheme of
$\SU_X$ defined above. Then \roster
\item The natural transformation in Theorem 1.2 induces a
transformation $$S\SU_X^{\natural}\to S\SU_X$$ such that $S\SU_X$
universally corepresents $S\SU_X^{\natural}$.
\item $S\SU_X$ is a projective, seminormal variety of dimension
$(r^2-1)(g-1)$. The number of irreducible components of $S\SU_X$
is the same with that of $\SU_X$, and its geometric points are
$$S\SU_X(\Bbb C)=\left\{\aligned &\text{$s$-equivalent classes of
semistable sheaves $\SE$}\\ &\text{of rank $r$ and degree $d$ with
$det(\SE|_{X^0})=\SL_0|_{X^0}$}
\endaligned\right\}.$$
\item There is an open subset $S\SU_X^s\subset S\SU_X$ whose
points represent isomorphism classes of stable sheaves with fixed
determinant $\SL_0$ on $X^0=X\setminus\{x_0\}.$ Locally in the
{\'e}tale topology on $S\SU_X^s$, there exists a universal sheaf
$\SE^{univ}$ such that if $\SE\in S\SU_X^{\natural}(S)$ whose
fibres $\SE_s$ are stable, then the pull-back of $\SE^{univ}$ via
$S\to S\SU^s_X$ is isomorphic to $\SE$ after tensoring with the
pull-back of a line bundle on $S$.\endroster\endproclaim

\demo{Proof} Firstly, it is easy to check (2). In fact, the
projectivity and the number of components follow from the
construction of $S\SU_X$. The semi-normality of $S\SU_X$ follows
from Corollary 1.5. To show that the set $S\SU_X(\Bbb C)$ consists
of the sheaves $F$ satisfying $det(F|_{X^0})=\SL_0|_{X^0}$, we
only need a fact that any line bundles $\SL_1$, $\SL_2$ on a
smooth projective curve $Y$ satisfy $\SL_1|_{Y\setminus\{y\}}
=\SL_2|_{Y\setminus\{y\}}$ if and only if
$\SL_1=\SL_2\otimes\SO_Y(ky)$ for some integer $k$.

To prove (1), let $G$ be the functor represented by $\SR^{ss}$ and
$F$ be the subfunctor defined by
$F(S)=\{(\SW_S\twoheadrightarrow\SE)\in
G(S)\,|\,det(\SE|_{X^0_S})=p_{X^0}^*\SL_0\}.$ Let
$Z\subset\SR^{ss}$ be the inverse image of $S\SU_X$ and $\Cal I_Z$
denote the ideal sheaf of $Z$. It is enough to show that for any
$(\SW_S\twoheadrightarrow\SE)\in F(S)$, the morphism
$\varphi_S:S\to \SR^{ss}$ factors through $\varphi_S:S\to Z\subset
\SR^{ss}.$  Namely, one has to prove $\varphi_S^*(\Cal I_Z)=0$,
where $\varphi_S^*:\SO_{\SR^{ss}}\to\varphi_{S*}\SO_S$. This is a
local problem, it is enough to show that for any $s\in S$ the
morphism $$\varphi_{S,s}:\hat S=Spec(\hat\SO_{S,s})\to \SR^{ss}$$
factors through $Spec(\hat\SO_{Z,\varphi_S(s)})$. By Proposition
1.3 and Corollary 1.5, there exists a complete noetherian local
$\Bbb C$-algebra $R$ and $\xi=\underarrow{\lim}(\xi_n),$ where
$$\xi_n=(\SW_{Spec(R/m^n)}\twoheadrightarrow\SE_n)\in
F(Spec(R/m^n)),$$ pro-represents the local deformation functor of
$F$ at $s$. Thus $$\varphi_{S,s}:\hat S=Spec(\hat\SO_{S,s})\to
\SR^{ss}$$ factors through $f:Spec(R)\to \SR^{ss}$, where the
pullback $f^*(\SW_{\SR^{ss}}\twoheadrightarrow\SE^{univ})$ of the
universal quotient is $\xi=\underarrow{\lim}(\xi_n)$, which is in
fact an element of $F(Spec(R))$ in our case (it is not true for
general theory). Thus $f=f_{\xi}:Spec(R)\to \SR^{ss}$ was given by
an element $\xi\in F(Spec(R))$. Now we can use Lemma 1.7 below,
which implies that $f_{\xi}^*(\Cal I_Z)=\sqrt{0}$, to prove that
$f_{\xi}$ (and thus $\varphi_{S,s}$) factors through $Z$ since $R$
is reduced.

Having shown (1) and (2), the proof of claim (3) is the same with
that of [Si].\enddemo

\proclaim{Lemma 1.7} Let $\xi\in F(S)$ and $f_{\xi}:S\to \SR^{ss}$
be the induced morphism by $\xi$, then $f_{\xi}(\wp)\in Z$ for any
point $\wp\in S$ (including non-closed points).\endproclaim

\demo{Proof} Let $\xi=(\SW_S\twoheadrightarrow \SF)\in F(S)$. When
$S$ is reduced, we have a stratification $S=\coprod S_a$ of $S$,
where  $S_a:=\{s\in S|\bold a(\SF_s)=a\}$ are locally closed
subschemes ($\bold a(\SF_s)$ was defined by
$\SF_s\otimes\hat\SO_{X,x_0}=\hat\SO_{X,x_0}^{\oplus\bold
a(\SF_s)}\oplus\hat m_{x_0}^{\oplus(r-\bold a(F_s))}$ ), then
using Lemma 2.7 of [S1] (for $S_a$) we get $f_{\xi}(S_a)\subset
Z.$ thus proves the lemma. For general $S$, we use the flattening
stratifications $S=\coprod S_i$ of Mumford (See Lecture 8 of [Mu])
such that the sheaves used in the construction of Lemma 2.7 of
[S1] are flat on $S_i$, then Lemma 2.7 goes through, thus the
lemma is proved for general $S$.\enddemo

\heading \S2 Degeneration of moduli space of semistable
$SL(r)$-bundles\endheading

Let $D$ be a complete discrete valuation ring with maximal idea
$m_D=(t)D$ and $\Cal X\to B=Spec(D)$ a flat family of proper
connected curves. Assume that the generic fibre $\Cal X_{\eta}$ is
smooth and the closed fibre $\Cal X_0$ is the curve $X$ discussed
in Section 1. Fix a relative ample line bundle $\SO_{\Cal X}(1)$
on $\Cal X$ such that $\SO_{\Cal X}(1)|_{\Cal X_0}=\SO_X(1)$. Let
$\SW_{\Cal X}=\SO_{\Cal X}(-N)^{\oplus P(N)}$ and
$Quot_{P(n)}(\SW_{\Cal X})\to B$ be the relative Grothendieck
quotient scheme. For any $B$-scheme $S$, we will write $\Cal X_S$
(resp. $\SW_{\Cal X_S}$) for $\Cal X\times_BS$ (resp. $\SO_{\Cal
X\times_BS}(-N)^{\oplus P(N)}$). Let $\SR^{ss}_{\Cal X}$ be the
open set of semistable sheaves whose quotient map induces
isomorphism $\SO_S^{P(N)}\cong H^0(\SE(N))$. It is well-known that
the relative moduli space  of semistable torsion free sheaves is
the relative good quotient $$\CD
 \SR^{ss}_{\Cal X}@>>> \SU_{\Cal X}:= \SR^{ss}_{\Cal X}//SL(P(N)) \\
 @VfVV               @V{\wt f}VV      @.    \\
  B          @=              B         @.
\endCD$$

Let $\SL$ be a line bundle on $\Cal X$ such that $\SL|_{\Cal
X_0}=\SL_0$ and $\SL_{\eta}=\SL|_{\Cal X_{\eta}}$. Let
$$\SR^{ss}_{\Cal X_{\eta}}(\SL_{\eta})\subset \SR^{ss}_{\Cal
X_{\eta}}$$ be the closed subscheme of sheaves with fixed
determinant $\SL_{\eta}$, and $$S\SU_{\Cal
X_{\eta}}:=\SR^{ss}_{\Cal X_{\eta}}(\SL_{\eta})//SL(P(N))$$ be the
moduli space of semistable bundles on $\Cal X_{\eta}$ with fixed
determinant $\SL_{\eta}.$  Let $$\SR^{ss}_{\Cal
X}(\SL):=\overline{\SR^{ss}_{\Cal X_{\eta}}(\SL_{\eta})} \subset
\SR^{ss}_{\Cal X},\quad S\SU_{\Cal X}:=\overline{S\SU_{\Cal
X_{\eta}}}\subset \SU_{\Cal X}\tag2.1$$ be the Zariski closure of
$\SR^{ss}_{\Cal X_{\eta}}(\SL_{\eta})$ and $S\SU_{\Cal X_{\eta}}$
in $\SR^{ss}_{\Cal X}$ and $\SU_{\Cal X}$ respectively. Then
$S\SU_{\Cal X}$ is the relative good quotient of $\SR^{ss}_{\Cal
X}(\SL)$. We will prove that $$\wt f:S\SU_{\Cal X}\to B$$
universally corepresents a moduli functor $S\SU_{\Cal
X}^{\natural}$, which is defined as $$S\SU_{\Cal X}^{\natural}(S)=
\left\{\aligned &\text{$\SO_S$-flat semistable sheaves $\SE$ on
$\Cal X_S:=\Cal X\times_BS$ of }\\ &\text{ rank $r$ and degree $d$
satisfying $det(\SE|_{\Cal X^0_S})=p_{\Cal X^0}^*\SL$}
\endaligned\right\},$$
where $\Cal X^0:=\Cal X\setminus\{x_0\}$ and $p_{\Cal X^0}:\Cal
X^0_S:=\Cal X^0\times_BS\to\Cal X^0$ denotes the projection.
Actually, we will prove that the functor $$\SR^{\natural}_{\Cal
X}(\SL)(S)=\left\{\aligned &\text{$\SO_S$-flat quotients
$(\SW_{\Cal X_S}\twoheadrightarrow \SE_S)$ }\\ &\text{ on $\Cal
X_S$ such that $\SE_S\in S\SU^{\natural}_{\Cal X}(S)$.}
\endaligned\right\}$$
is represented by the closed subscheme $\SR^{ss}_{\Cal
X}(\SL)\subset\SR^{ss}_{\Cal X}.$

\proclaim{Lemma 2.1} The functor $\SR^{\natural}_{\Cal X}(\SL)$ is
represented by the closed subscheme $\SR^{ss}_{\Cal
X}(\SL)\subset\SR^{ss}_{\Cal X}$ and the restriction of universal
quotient on $\Cal X\times_B\SR^{ss}_{\Cal X}$.\endproclaim

\demo{Proof} For any $(\SW_{\Cal X_S}\twoheadrightarrow
\SE_S)\in\SR^{\natural}_{\Cal X}(\SL)(S)$, there is a unique
morphism $$\varphi_S:S\to\SR^{ss}_{\Cal X}$$ such that pullback of
the universal quotient $(\SW_{\Cal X_{\SR^{ss}_{\Cal
X}}}\twoheadrightarrow \SE^{univ}_{\SR^{ss}_{\Cal X}})$ is
$(\SW_{\Cal X_S}\twoheadrightarrow \SE_S)$. It is enough to show
that $\varphi_S$ is factorized through $$\CD
S@>\varphi_S>>\SR^{ss}_{\Cal X}(\SL)@>>>\SR^{ss}_{\Cal X}\\ @VVV
@Vf_{\SL}VV @VfVV  \\ B@= B@= B \endCD$$ This is true at the
generic fibre, we only need to check the case when $S$ is defined
over $0\in B.$ In the proof of Theorem 1.6 (1), we have show that
$\varphi_S$ is factorized through (note that we are at the case of
$\Cal X\times_BS=\Cal X_0\times S$) $$
Z=\{(\SW\twoheadrightarrow\SE)\in \SR^{ss}_{\Cal X_0}\,|\,
det(\SE|_{X^0})=p^*_{X^0}(\SL_0)\}.$$ Thus the lemma follows the
equality: $f_{\SL}^{-1}(0):=\SR^{ss}_{\Cal X}(\SL)_0=Z,$ which we
will prove in the next lemma.\enddemo

\proclaim{Lemma 2.2} $f_{\SL}^{-1}(0):=\SR^{ss}_{\Cal
X}(\SL)_0=Z$.\endproclaim

\demo{Proof} We first prove that $\SR^{ss}_{\Cal X}(\SL)_0\subset
Z$. Let $\xi_0=(\SW\twoheadrightarrow\SE_0)\in\SR^{ss}_{\Cal
X}(\SL)_0$, we need to show that
$det(\SE_0|_{X^0})=p^*_{X^0}(\SL_0)$.

By the definition of $\SR^{ss}_{\Cal X}(\SL)$, there exists a
complete discrete valuation ring $\wt D$ dominant $D$, and a
morphism $T:=Spec(\wt D)@>\varphi>>\SR^{ss}_{\Cal X}(\SL)$ over
$B$ such that $$\varphi(\wt\eta)\in\SR^{ss}_{\Cal
X_{\eta}}(\SL_{\eta}),\quad\varphi(\wt 0)=\xi_0.$$  Namely, there
is a $\SO_T$-flat torsion free sheaf $\SE$ on $\Cal X\times_BT$
such that $det(\SE|_{\Cal X_{\eta}})=\SL_{\eta}$ and $\SE|_{\Cal
X_0}=\SE_0.$

Let $\varpi:\wt{\Cal X}_T \to \Cal X_T$ be the desingularization
of $\Cal X_T$ at $x_0$. The exceptional divisor
$\varpi^{-1}(x_0)=\sum E_i$ is a chain of $(-2)$-curves, and the
special fibre $\wt{\Cal X}_0$ of $\wt{\Cal X}_T\to T$ at $\wt 0$
is $X_1+X_2+\sum E_i$, and $X_1\cap (X_2+\sum E_i)=\{x_1\}$,
$X_2\cap (X_1+\sum E_i)=\{x_2\}.$ Identifying $\wt{\Cal
X}_T\setminus\varpi^{-1}(x_0)\cong\Cal X_T\setminus\{x_0\}=\Cal
X^0_T$, we can extend $det(\SE|_{\Cal X^0_T})$ into a line bundle
$\overline{det(\SE|_{\Cal X^0_T})}$ on $\wt{\Cal X}_T$, which
satisfies that $(\varpi^*p^*_{\Cal X}\SL)|_{\wt{\Cal X}_{\eta}}=
\overline{det(\SE|_{\Cal X^0_T})}|_{\wt{\Cal X}_{\eta}}.$ Thus one
has $$\overline{det(\SE|_{\Cal X^0_T})}=(\varpi^*p^*_{\Cal
X}\SL)\otimes\SO_{\wt{\Cal X}_T}(V),$$ where $V\subsetneqq\wt{\Cal
X}_0$ is a vertical divisor. Since $\SO_{\wt{\Cal X}_T}(\wt{\Cal
X}_0)$ is trivial, there are integers $n_1$, $n_2$ such that
$\SO_{\wt{\Cal X}_T}(V)|_{X_1}=\SO_{X_1}(n_1x_1)$, $\SO_{\wt{\Cal
X}_T}(V)|_{X_2}=\SO_{X_2}(n_2x_2).$ Hence
$$det(\SE_0|_{X^0})=det(\SE|_{\Cal X^0_T})|_{X^0}=
\overline{det(\SE|_{\Cal X^0_T})}|_{X^0}=p^*_{X^0}(\SL_0).$$

We are left to prove $Z\subset\SR^{ss}_{\Cal X}(\SL)_0$. Since the
locus $\SR_0\subset f^{-1}(0)=\SR^{ss}_{\Cal X_0}$ of locally free
sheaves is a dense open subset of $\SR^{ss}_{\Cal X_0}$, we only
need to check that $$\SR_0\cap Z\subset\SR^{ss}_{\Cal X}(\SL)_0$$
Let $\xi_0=(\SW\twoheadrightarrow\SE_0)\in\SR_0\cap Z$, then
$\SE_0$ is a vector bundle on $\Cal X_0=X$ and
$det(\SE_0|_{X^0})=p^*_{X^0}(\SL_0)$, which implies that (for
$i=1,2$)
$$det(\SE_0)|_{X_i}=det(\SE_0|_{X_i})=\SL_0|_{X_i}\otimes\SO_{X_i}(n_ix_0),\quad
n_1+n_2=0.$$ Let $\wt\SL=\SL\otimes\SO_{\Cal X}(X_1+(n_1+1)X_2)$
and $\wt\SL_0=\wt\SL|_{\Cal X_0}$. Then, since $\SO_{\Cal
X}(X_1+X_2)\cong\SO_{\Cal X}$, we have $\SO_{\Cal
X}(X_1+(n_1+1)X_2)|_{X_i}=\SO_{X_i}(n_ix_0)$ and
$\wt\SL_0|_{X_i}=det(\SE_0)|_{X_i},$ which implies that
$$det(\SE_0)=\wt\SL_0=\wt\SL|_{\Cal X_0}.$$

Let $d_i=deg(\SE_0|_{X_i})$ and $J_{X_i}^{d_i}$ ($i=1,2$) be the
Jacobian variety of line bundles of degree $d_i$ on $X_i$. Let
$J(\Cal X_{\eta})$ be the Jacobian varirty of line bundles of
degree $d$ on $\Cal X_{\eta}$. Then $J(\Cal X_{\eta})$ can be
compactified into a relative Jacobian variety $$J(\Cal X)\to B$$
such that its special fibre $J(\Cal X)_0=J^{d_1}_{X_1}\times
J^{d_2}_{X_2}$ (Otherwise, we can modify $J(\Cal X)$ through an
isomorphism by tensoring a suitable line bundle $\SO_{\Cal
X}(k_1X_1+k_2X_2)$ on $\Cal X$). Let $\SR^0_{\Cal
X}\subset\SR^{ss}_{\Cal X}$ be the locus of locally free sheaves.
By taking determinant, we have a morphism $$Det:\SR^0_{\Cal
X}\setminus\bigcup_{\chi_i\neq
d_i+r(1-g_i)}\SR_0^{\chi_1,\chi_2}\to J(\Cal X).$$ The line bundle
$\wt\SL$ on $\Cal X$ gives a section $\sigma:B\to J(\Cal X)$ such
that $\sigma(0)=\sigma(B)\cap J(\Cal X)_0$. Let
$\SR^{\sigma}:=Det^{-1}(\sigma(B))\to B$ and
$\SE^{\sigma}=\SE^{univ}|_{\Cal X\times_B\SR^{\sigma}}.$ Then
$$det(\SE^{\sigma})=p^*_{\Cal X}(\wt\SL)\otimes
p^*_{\SR^{\sigma}}(\SK)$$ for some line bundle $\SK$ on
$\SR^{\sigma},$ and $\xi_0\in\SR^{\sigma}$. Thus, for any open set
$\xi_0\in U$ such that $\SK$ is trivial on $U\to B$, we have
$U\cap\SR^{ss}_{\Cal X}(\SL_{\eta})\neq\emptyset.$ Thus
$$\xi_0\in\overline{\SR^{ss}{\Cal X}(\SL_{\eta})},$$ which proves
the lemma.
\enddemo

Let $\bl_{\Cal X}$ be the category of Artinian local $D-$algebras.
For any point $$e=(\SW\twoheadrightarrow \SE_0)\in\SR^{ss}_{\Cal
X}(\SL)_0,$$ let $F_{\Cal X}:\bl_{\Cal X}\to \bset$ denote the
local deformation functor of $\SR^{ss}_{\Cal X}(\SL)$ at $e$.
Similarly, we define a functor $T_{\Cal X}:\bl\to \bset$ by
$$T_{\Cal X}(A)=\{\text{Isomorphism classes of $A$-flat torsion
free $(\SO_{X,x_0}\otimes_D A)$-modules}\}$$ such that $T(\Bbb
C)=\{\SE_0\otimes\SO_{X,x_0}\}.$ There is a morphism of functors
$$\phi: F_{\Cal X}\to T_{\Cal X}$$ defined by $\phi((\SW_{\Cal
X}\twoheadrightarrow \SE_{\Cal X_ A}))=\SE_{\Cal X_
A}\otimes(\SO_{\Cal X,x_0}\otimes_D A).$ Similar with Lemma 1.4,
we have

\proclaim{Lemma 2.3} The morphism $\phi: F_{\Cal X}\to T_{\Cal X}$
is formally smooth.\endproclaim

\demo{Proof} The same arguments with Lemma 1.4, we just remark an
easy fact: Let $B$ be a flat $D$-algebra and $M$ be a $B$-module,
flat over $D$. Assume that $M\otimes_DD/m$ is a free
$B/mB$-module. Then $M$ is a free $B$-module.\enddemo

Let $\bold X=(x_{ij})_{(r-a)\times(r-a)}$  and $\bold
Y=(y_{ij})_{(r-a)\times(r-a)}$ be the $(r-a)\times(r-a)$ matrics,
and $$\Cal Z=\text{Spec$\frac{D[\bold X,\bold Y]}{(\bold
X\cdot\bold Y-t,\bold Y\cdot\bold X-t)}$}.$$ Let $0\in\Cal Z$ be
the point defined by the idea $(\bold X,\bold Y, t)\SO_{\Cal Z}.$
Then we have

\proclaim{Lemma 2.4} Let $e=(\SW\twoheadrightarrow
\SE_0)\in\SR^{ss}_{\Cal X}(\SL)_0$ such that $\bold a(\SE_0)=a$.
Then there exist $\ell_1$ and $\ell_2$ such that
$$\hat\SO_{\SR^{ss}_{\Cal X}(\SL),e}[[u_1,\cdots,u_{\ell_1}]]
\cong\hat\SO_{\Cal Z,0}[[v_1,\cdots,v_{\ell_2}]].$$ In particular,
when $\Cal X$ is regular, $\SR^{ss}_{\Cal X}(\SL)$ is
regular.\endproclaim

\demo{Proof} This is the consequence of Lemma 2.3 and the results
of [F2] and [NS].\enddemo

\proclaim{Theorem 2.5} Let $\Cal X\to B$ be a regular scheme with
closed fibre $\Cal X_0=X$ and a fixed relative ample line bundle
$\SO_{\Cal X}(1)$. Let $\SL$ be a line bundle on $\Cal X$ such
that $\SL|_{\Cal X_0}=\SL_0$ and let $\wt f: S\SU_{\Cal X}\to B$
be defined as $(2.1)$. Then \roster
\item There is a natural transformation $S\SU_{\Cal X}^{\natural}\to S\SU_{\Cal X}$
such that $\wt f:S\SU_{\Cal X}\to B$ universally corepresents
$S\SU_{\Cal X}^{\natural}$.
\item $\wt f:S\SU_{\Cal X}\to B$ is a flat family of projective varieties of
dimension $(r^2-1)(g-1)$, whose general fibres $\wt f^{-1}(b)$ are
moduli spaces $S\SU_{\Cal X_b})$ of semistable vector bundles of
rank $r$ and degree $d$ with fixed determinant $\SL_b=\SL|_{\Cal
X_b}$, and its closed fibre $\wt f^{-1}(0)$ is the variety
$S\SU_X$ in Theorem 1.6.
\item The scheme $S\SU_{\Cal X}$ is normal with only rational singularities,
and the locus $S\SU_{\Cal X}^s$ of stable sheaves is smooth.
\endroster\endproclaim

\demo{Proof} (1) follows from Lemma 2.1 and the same arguments of
[Si]. (2) is clear. (3) follows from Lemma 2.4 and general
theorems in GIT.\enddemo

\enddocument